\documentclass[reqno]{amsart}
\usepackage{amssymb}
\usepackage{graphicx}
\usepackage[left]{lineno}
\usepackage{placeins}
\usepackage{float}
\usepackage{xcolor}

\newtheorem{theorem}{Theorem}[section]

\newtheorem{proposition}[theorem]{Proposition}
\newtheorem{corollary}[theorem]{Corollary}

\theoremstyle{definition}
\newtheorem{definition}[theorem]{Definition}

\theoremstyle{remark}

\numberwithin{equation}{section}

\usepackage[margin=2.2cm]{geometry}
\DeclareMathOperator*{\argmin}{arg\,min}
\DeclareMathOperator*{\dom}{dom}

\DeclareMathOperator*{\Lip}{Lip}

\title{On the geometry of sharp minima}

	\author[]{Alberto Dom\'inguez Corella}
    
	\address{Institut f\"{u}r Stochastik und Wirtschaftsmathematik, VADOR E105-04, TU Wien, Wiedner Hauptstra{\ss}e 8, A-1040 Wien, \"Osterreich}

    \address{Institut f\"{u}r Mathematik und Wissenschaftliches Rechnen, Universit\"{a}t Graz, Heinrichstra{\ss}e 36, A-8010 Graz, \"Osterreich\\}

	\email{alberto.of.sonora@gmail.com}
	
\thanks{The author was supported by the Austrian Science Foundation (FWF) under grants P 36344-N and F 100800.}

\begin{document}

\begin{abstract}
   We give characterizations of sharp minimizers that emphasize their geometric properties. These include tilt invariance and weak upper gradient conditions. We relate sharp minimality to cusps in nonsmooth manifolds when interpreted locally as graphs and connect the rolling of tangent planes to the tilting of functions.
\end{abstract}

\maketitle

\section{Introduction}
	Sharp minimizers are a topic of interest in optimization, e.g., in convex programming \cite{BF_1993} and  convergence analysis \cite{BF_1995}. Their relevance lies in their robustness to perturbations and in enabling finite-time convergence for certain algorithms. In words, a minimizer is sharp if the function increases at least linearly away from it, ruling out any flattening around it. From the geometric point of view, this corresponds to the minimizer being the vertex of a cone that lies below the function.
	\smallbreak \noindent
	The manuscript is devoted to a geometric interpretation of sharp minimizers, backed by mathematical statements. The results comprise characterizations that emphasize geometric features arising from the rolling of tangent planes. Below, here in the introduction, we provide a preamble to the main results by discussing two-dimensional nonsmooth manifolds and their relation to sharp minimality when interpreted locally as graphs. We then comment on the related literature in the context of our contributions. The rigorous mathematical statements are presented in the next section, and the proofs are deferred to the final section.
	\subsection{Cones and tilting}
	Let us begin our journey by talking about cones. What makes cones so special? Why are they different from, say, a sphere? Well, every point is indistinguishable on a sphere, whereas a cone has one very special point (its vertex).
	\vspace{-0.6cm}
	\begin{figure}[H]
\includegraphics[width=0.75\textwidth]{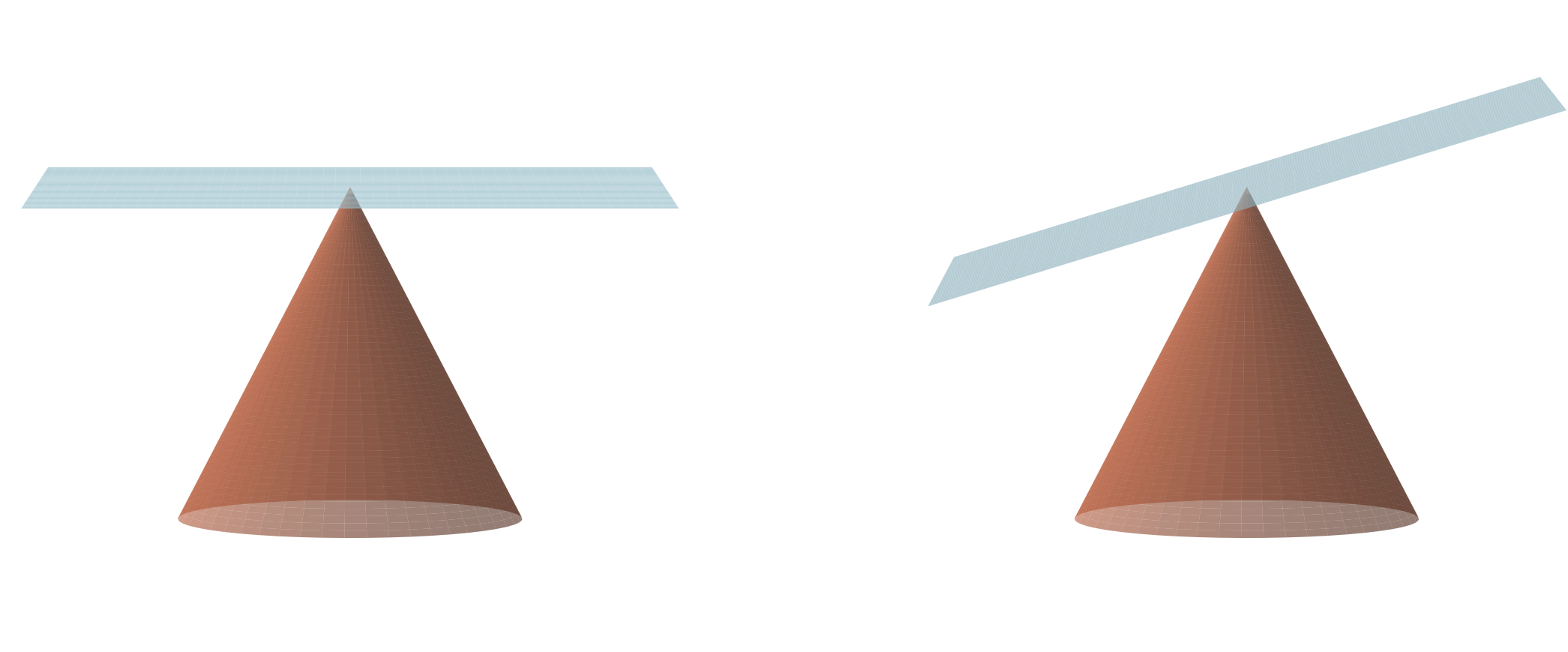}
    \vspace{-0.7cm}
		\caption{Rolling of the tangent plane}
		\label{fig:tangents}
	\end{figure}\noindent
	To start with, we can think of a cone as a 2-dimensional manifold. At its vertex, there are several tangent planes, including one that is perpendicular to the cone. If we slightly rolled this tangent plane, it would continue to act as a tangent plane;  this is unlike points elsewhere on the cone, where rolling the tangent plane would not maintain tangency. A natural question arises: if a point on a two-dimensional surface has the property that one of its tangent planes can be slightly perturbed in every direction while still remaining tangent at that point, does this imply that the surface is locally a cone with the point as its vertex? While the answer is generally no, we will see that  such a property does imply that the surface can be locally enclosed within the convex hull of a cone.
	\smallbreak 
	\noindent
	We can think of a cone as the graph of a function
	\( c: \mathbb{R}^2 \to \mathbb{R} \) of the form
	\begin{align}\label{coneform}
		c(x) =\cot(\alpha/2) \| x  - p \|  + \big(\tan(\beta) - \cot(\alpha/2)\big) \langle v, x - p \rangle,
	\end{align}
	where \( \alpha\in(0,\pi) \) is the maximum angle between two generatrix lines, the vertex is located at $p\in\mathbb R^2$, \(v\in S^1\) is the direction vector where the cone is leaned, and $\beta\in(0,\pi/2)$ is the angle between the cone and the $xy$-plane.  After applying suitable translations and rotations, we may reduce cones to the simpler form \( c(x) = \gamma \| x \| \) for some \( \gamma > 0 \); in this form, the vertex is at the origin, and the perpendicular tangent plane is the $xy$-plane.  
	\smallbreak 
	\noindent
	Given a 2-dimensional manifold admitting a tangent plane at a point, we locally represent it as the graph of a proper lower semicontinuous function 
	$f: \mathbb{R}^2 \to \mathbb{R} \cup \{+\infty\}$ 
	that takes the value $+\infty$ outside a closed neighborhood and has the origin as unique minimizer; here, $(0, f(0))\in\mathbb R^3$ represents the point admitting the tangent plane.  
	An alternative way of thinking of the rolling of the tangent plane at the point is  by fixing the plane and tilting the graph itself. In this interpretation, tilting the graph in a direction \( v\in\mathbb R^2 \) corresponds to subtracting the term \( \langle v, \cdot \rangle \) from the function \( f \). 
	\smallbreak\noindent
	If we slightly tilt a cone \( c(x) = \gamma \|x\| \), the \( xy \)-plane remains a tangent plane to the perturbed cone. Now, suppose that when we tilt the function \( f \) in small directions \( v \in \mathbb{R}^2 \), the $xy$-plane remains below the function, i.e.,  the origin remains the unique minimizer; in symbols, 
	\[
	\argmin_{y\in \mathbb{R}^2} \{ f(y) - \langle v,y \rangle \} = \{0\}.
	\]  
	This corresponds to saying that at the given point on the surface, the tangent plane can be slightly rolled while still remaining tangent. Theorem \ref{Thm1} in the next section states that this is equivalent to the existence of a cone \( c \) such that \( f \geq c \). Intuitively, this means that near the point admitting the tangent plane, the surface lies inside the convex hull of a cone; this is because we view the graph of $f$ as a neighborhood of the surface, and the conclusion $f \ge c$ says exactly that the epigraph of $f$ is contained in that of $c$.
\vspace{-1cm}
	\begin{figure}[H]
		\includegraphics[width=0.70\textwidth]{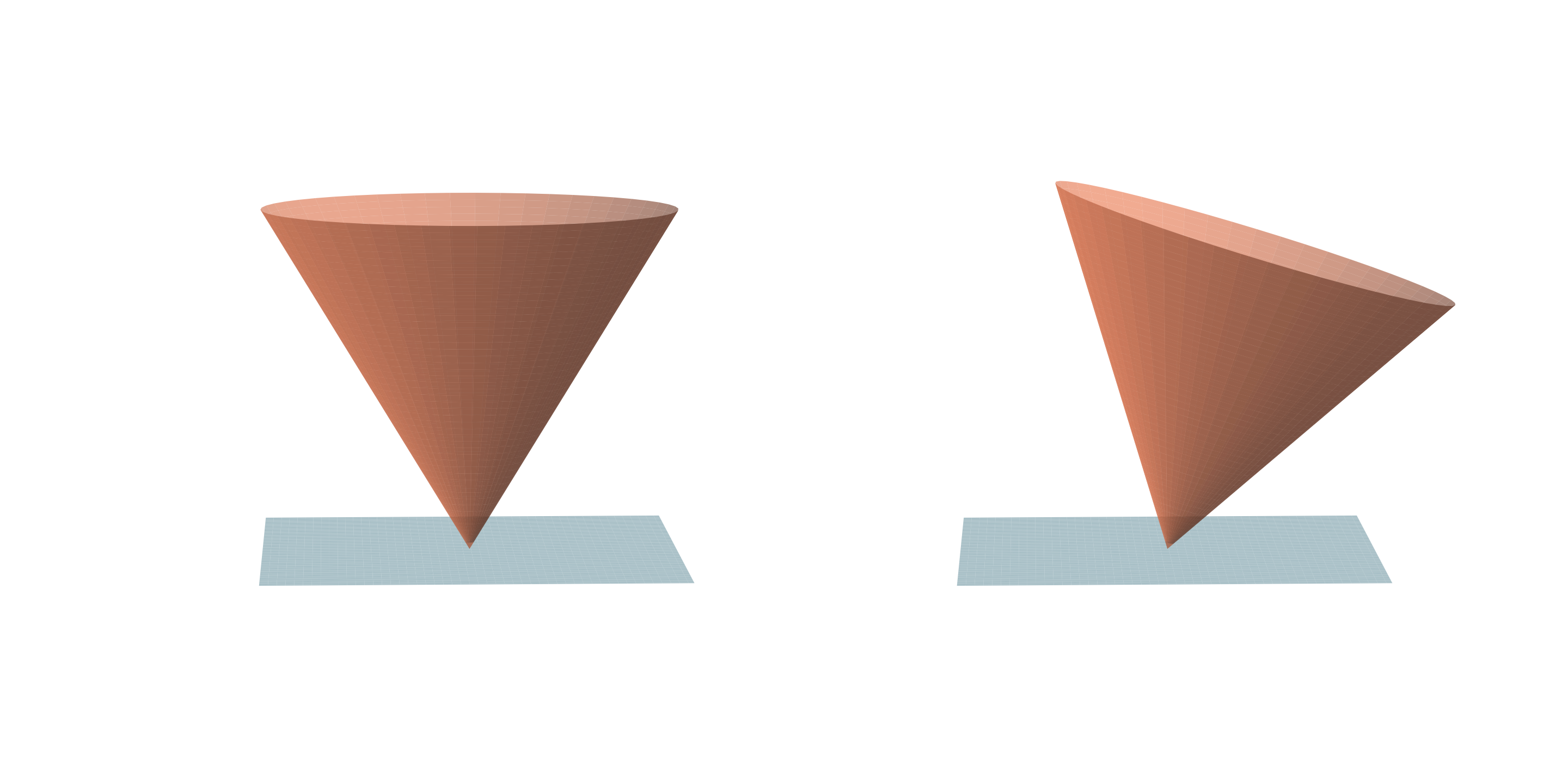}
        \vspace{-1.4cm}
		\caption{Tilting of the cone}
		\label{fig:cones}
	\end{figure}
    \vspace{-0.3cm}\noindent 
In relation to cusps, Proposition \ref{P1} gives a characterization of sharp minimality based on a weak upper gradient (the global metric slope). This result leads to the geometric conclusion that the  graph of $f$ does not flatten around $(0,f(0))$ if and only if it is contained within the convex hull of a cone. Lastly, Proposition \ref{pappendix} shows that one can convexify a local model (the epigraph of $f$) of the surface around a point and still detect the invariance of the tangent plane under small tilts at that point.
	\subsection{On the literature}
	The concept of a sharp minimizer was introduced by Polyak in  \cite[p. 205]{P_1987}.  This notion was extended by Ferris in \cite{F_1988} to encompass cases with multiple minimizers, referred to as weak sharp minimizers. Subsequently, in \cite{BF_1993}, Burke and Ferris demonstrated the significance of weak sharp minimizers in convex programming and convergence analysis. The primary motivation for studying this concept lies in its impact on sensitivity analysis and the convergence behavior of a wide range of optimization algorithms. For a comprehensive account and detailed references, we refer to \cite{BD_2002,BD_2005,BD_2009}. Over the years, several works have provided sufficient conditions and characterizations of sharp minimizers. Notable contributions include \cite{A_2006, KZ_2003, U_2014, U_2015, W_1994, ZMX_2012, ZW_2012, ZX_2011}. It is worth noting that the concept of tilt-stability (Lipschitz stability of minimizers in a neighborhood) given in \cite[Definition 1.1]{PR_1998} is different from the tilting invariance considered in this paper (which can be regarded as a  much stronger form of stability).
	\smallbreak
	\noindent
Our main contribution is the geometric interpretation given in the previous subsection, along with the tilting results stated in Theorem \ref{Thm1} and Corollary \ref{Cor1} in Section \ref{Sec_gbl}, which support this interpretation in normed spaces. A secondary contribution is the extension of these results to metric spaces and to sharp local   minimizers in Subsection \ref{Sec_loc}.
Other minor contributions are Propositions \ref{P1} and \ref{P2}, which have appeared several times in the literature before, in other forms, contexts, or under different assumptions.
	\smallbreak 
	\noindent
	In \cite[Section 7.1.3]{P_1987}, the concept of superstability was introduced, and it was shown that sharp minimizers are superstable. There, superstability refers to the preservation of minimizers under a certain class of small convex perturbations. In Corollary \ref{Cor2}, we provide more comprehensive results on convex perturbations. 
	In \cite[Proposition 2]{U_2015}, a different interpretation of superstability is considered, focusing on local minimizers and general perturbations. The generality of these perturbations allows for a characterization of sharp local    minimality by exploiting the many ways a function can be perturbed. In contrast, our focus is on the preservation of global minimizers, even when addressing local conditions—this is essential for the geometric interpretation introduced earlier.  We also allow weaker perturbation families to characterize sharp minimizers, e.g., Theorem \ref{Thm1} shows  sharp minimality can be characterized via the action of the topological dual on the original functional.
	\smallbreak 
	\noindent
	In Propositions \ref{P1} and \ref{P2}, we give characterizations of sharp minimality based on metric slopes. The statements are not entirely new, but they do not seem to appear elsewhere in the literature in the form in which we present them. In this direction, there have been several related results concerning subdifferentials; see \cite[Theorem 5.2]{SW_1999} for Euclidean spaces and \cite[Theorem 2.2]{KZ_2003} for Asplund spaces. See also \cite[Theorem 2.4]{U_2014} and \cite[Proposition 1]{U_2015} for further characterizations. We note that the condition that the slope is bounded below by a positive constant appears in the study of Eikonal equations in metric spaces, where it corresponds to the notion of Monge super-solutions; see, for example, \cite{GS_2015, GHN_2015, LSZ_2021}. In that setting, such conditions are used to analyze viscosity or metric solutions to Hamilton–Jacobi equations; see also \cite{DLS_2024, DS_2022, LT_2024}. Finally, a complementary geometric approach to optimization, involving Morse-theoretic ideas as well as parametric and stability analysis, was developed in the finitely constrained setting in \cite{JJT_1983,JW_1990,JW_1991}. Related developments for semi-infinite optimization and optimal control can be found in \cite{JTW_1992,W_1998}.

    \section{Geometry of sharp minimizers}
	The main results are given in this section and the  proofs are deferred to the next one. 
	In the first part, we give   characterizations for global sharp minimizers of functions over normed spaces. 
	In the second part, we consider sharp local    minimizers of functionals defined on metric spaces.
	\subsection{Sharp global minimizers in normed spaces}\label{Sec_gbl}
	Throughout this subsection, we consider a real normed space \((X, \|\cdot\|)\) and a proper function \(f:X\to \mathbb{R}\cup\{+\infty\}\).
	\smallbreak 
	We fix a reference point $\bar x\in\dom f$. 
	\begin{definition}
		We say  $\bar x$ is a  sharp minimizer if there exists $\gamma>0$ such that 
		\begin{align}\label{growth}
			f(x)\ge f(\bar x) +\gamma \|x-\bar x\|\quad \forall x\in X.
		\end{align}
		A number $\gamma>0$ satisfying  (\ref{growth}) is called a sharpness modulus. 
	\end{definition}
	\noindent 
	From a geometric point of view, a function with a sharp minimizer lies above a cone. Intuitively, if a function is bounded below by such a cone, slight tilting in any direction will not shift the minimizer. This also begs the question: if slight enough tilting does not change the minimizer, can a cone be fit below the function? The answer is positive and  condensed in the following theorem.
	\begin{theorem}\label{Thm1}
		Let $\gamma>0$ be given. The following statements are equivalent.
		\begin{itemize}
			\item[$(i)$] For any $\xi\in X^*$, 
			\[
			\|\xi\|<\gamma\quad\implies\quad\argmin_{y\in X}\big\{f(y)-\langle \xi,y\rangle\big\}=\{\bar x\}.
			\]
			
			\item[$(ii)$] There exists a dense set $\Theta\subseteq X^*$ such that, for any $\xi\in\Theta$, 
			\begin{align*}
				\|\xi\|<\gamma\quad\implies\quad\bar x\in \argmin_{y\in X}\big\{f(y) -\langle \xi, y\rangle\big\}.
			\end{align*}
			
			\item[$(iii)$] The point $\bar x$ is a sharp minimizer with modulus $\gamma$.
			
		\end{itemize}
	\end{theorem}
	\noindent 
	This result shows that sharp minimality can be determined by how the dual of the underlying space perturbs the original function. The proof of the previous result employs the Hahn-Banach theorem; it is used to construct adequate objects in the dual of the underlying space to perturb the original function.
	\smallbreak 
	\noindent 
    The  tilt invariance of minimizers can also be interpreted as stability with respect to linear perturbations; then a natural question is: for what other perturbations besides linear ones does this stability hold? It turns out that one can consider general Lipschitz perturbations.
	\begin{corollary}\label{Cor1}
		Let $\gamma>0$ be given. The  following statements are equivalent.
		\begin{itemize}
			\item[$(i)$] For any $\xi\in X^*$, 
			\[
			\|\xi\|<\gamma\quad\implies\quad\argmin_{y\in X}\big\{f(y)-\langle \xi,y\rangle\big\}=\{\bar x\}.
			\]
			\item[$(ii)$] For any Lipschitz function $\zeta: X\to\mathbb R$, 
			\[
			\Lip \zeta<\gamma\quad\implies\quad\argmin_{y\in X}\big\{f(y)+\zeta(y)\big\}=\{\bar x\}.
			\]
		\end{itemize}
	\end{corollary}
	\noindent 
	We  now give a characterization of sharp minimality in terms of weak upper gradients. 
	The global slope of $f$ at $x\in X$ is given by 
	\begin{align*}
		|\partial f|(x):= \sup_{y\in X\setminus\{x\}}\frac{\max\big\{f(x) - f(y),0\big\}}{\|x-y\|}.
	\end{align*}
	This notion can be used as a measure of minimality since  $|\partial f|(\bar x)=0$ if and only if $\bar x$ is a minimizer of $f$. If $f$ is convex, the global slope at a point coincides with the minimal norm subgradient of the function at that point. 
    \smallbreak\noindent
	The next proposition shows that  a point is a sharp minimizer  if and only if every other point has a global slope bounded away from zero. Intuitively, this means the function does not flatten around the minimizer.
	\begin{proposition}\label{P1}
		Suppose $X$ is a Banach space and  $f:X\to\mathbb R\cup\{+\infty\}$ is lower semicontinuous and bounded from below. 
		Let $\gamma>0$ be given. Then $\bar x$ is a sharp minimizer  with modulus $\gamma$ if and only if 
		\begin{align*}
			|\partial f|(x)\ge\gamma \quad  \forall x\in X\setminus\{\bar x\}. 
		\end{align*}
	\end{proposition}\noindent
	Observe that lower semicontinuity is necessary in the previous proposition. For example, consider the function  \( f(x) := 1 - |1-\|x-\bar x\|| \) if \(\|x-\bar x\|<2\) and \(+\infty\) otherwise. A similar example can show the need of completeness. 
    \smallbreak\noindent
   Since the notion of slope can be stated in purely metric settings \cite[Subsection 1.2]{AGS_2008}, it is natural to expect Proposition \ref{P1} to extend to complete metric spaces; this is indeed the case for global sharp minimizers, since the proof in Subsection \ref{pss} relies only on Ekeland's variational principle.
   \smallbreak\noindent
   The Legendre--Fenchel transform of $f$ is the function $f^*:X^*\to \mathbb R\cup\{+\infty\}$ given by $f^*(\xi) := \sup_{x \in X} \bigl\{ \langle \xi, x \rangle - f(x) \bigr\}$. 
The double Legendre--Fenchel transform, or biconjugate, of $f$ is the function $f^{**}:X\to \mathbb R\cup\{+\infty\}$ defined by
\[
f^{**}(x) := \sup_{\xi \in X^*} \bigl\{ \langle \xi, x \rangle - f^*(\xi) \bigr\}.
\]
Note that \(f^{**}:X\to \mathbb R\cup\{+\infty\}\)  coincides with the lower semicontinuous convex envelope of \(f\), i.e., the largest lower semicontinuous convex function \(g:X\to\mathbb{R}\cup\{+\infty\}\) such that \(g\le f\).
\smallbreak\noindent
The following proposition shows that sharp minimality is preserved under biconjugation.
\begin{proposition}\label{pappendix}
Suppose $\bar x \in \argmin_{y\in X} f(y)$ and let  $\gamma > 0$. Then $\bar x$ is a sharp minimizer of $f : X \to \mathbb{R} \cup \{+\infty\}$ with modulus $\gamma$ if and only if it is a sharp minimizer of $f^{**} : X \to \mathbb{R} \cup \{+\infty\}$ with modulus $\gamma$.
\end{proposition}
	\subsection{Sharp local  minimizers in metric spaces}\label{Sec_loc}
	Throughout this subsection, we consider a complete metric space \((\mathcal M,d)\) and a proper lower semicontinuous function \(\mathcal J:\mathcal M\to\mathbb{R}\cup\{+\infty\}\) bounded from below.
	\smallbreak
	We fix a reference minimizer $\bar u\in \mathcal M$ of $\mathcal J$.
	\begin{definition}
		We say $\bar u$ is a \textit{sharp local    minimizer}   if there exist $\delta\in(0,+\infty]$ and $\gamma>0$ such that 
		\begin{align}\label{slmy}
			\mathcal J(u)\ge \mathcal J(\bar u) +\gamma\hspace*{0.02cm}d(u,\bar u)\quad \forall u\in\mathbb B\big(\bar u,\delta\big).
		\end{align}
		We say that $\delta\in(0,+\infty]$ and $\gamma>0$ are  parameters of sharp minimality.
	\end{definition}
\noindent 
	It was seen in the previous subsection that sharp minimality can be characterized by the way Lipschitz mappings perturb the function in normed spaces. This is also the case in complete metric spaces. 
	\begin{theorem}\label{Thm2}
		Let $\gamma>0$ be given. The following statements are equivalent. 
        \begin{itemize}
            \item[$(i)$] There exists $\delta\in(0,+\infty]$ such that, for every Lipschitz function $\zeta:\mathcal M\to\mathbb R$,
		\begin{align*}
			\Lip\zeta<\gamma\quad\implies\quad \argmin_{v\in\mathcal M}\big\{\mathcal J(v) + \zeta(v)\big\}\cap\mathbb B(\bar u,\delta)\subseteq \{\bar u\}.
		\end{align*}

        \item[$(ii)$]  There exists $\delta\in(0,+\infty]$ such that  $\bar u$ is a sharp local    minimizer with parameters $\delta$ and $\gamma$
        \end{itemize}
	\end{theorem}
	\noindent 
    The proof of the preceding theorem relies on Ekeland’s variational principle. We note that, in the previous section, the analogous result for sharp global minimizers followed directly from the tilt-invariance characterization in Theorem \ref{Thm1} (itself based on the Hahn–Banach theorem). In the metric-space setting, however, no such direct route is available, and a different argument is required.
    \smallbreak\noindent
	The local metric slope of a proper  function $\varphi:\mathcal M\to\mathbb R\cup\{+\infty\}$ at a point $u\in\mathcal M$ is given by 
	\begin{align*}
		|\nabla \varphi|(u)
:=\begin{cases}
\displaystyle\limsup_{v\to u}\displaystyle\frac{\max\{\varphi(u)-\varphi(v),0\}}{d(u,v)}
& \text{if \(u\) is not isolated} \\
0 & \text{if \(u\) is isolated}.
\end{cases}
	\end{align*}
	The metric slope can be used as a measure of \textit{criticality}; a point for which the metric slope is zero is said to be critical. One can easily prove that the metric slope must vanish at  local minimizers. 
    \smallbreak\noindent 
   Next we show that, in spaces with nonpositive curvature, sharp minimality can be characterized in terms of how geodesically convex functionals perturb the original functional (as measured by the local slope). Preliminaries on nonpositively curved spaces and geodesic convexity are given in the appendix.
	\begin{corollary}\label{Cor2}
		Let $\gamma>0$ be given. If $(\mathcal M, d)$ is an Hadamard space, then the following statements are equivalent.
		\begin{itemize}
			\item[$(i)$] There exists $\delta\in(0,+\infty]$ such that, for every Lipschitz function $\zeta:\mathcal M\to\mathbb R$,
			\begin{align*}
				\Lip\zeta<\gamma\quad\implies\quad \argmin_{v\in\mathcal M}\big\{\mathcal J(v) + \zeta(v)\big\}\cap\mathbb B(\bar u,\delta)\subseteq \{\bar u\}.
			\end{align*}
			
			\item[$(ii)$] There exists $\delta\in(0,+\infty]$ such that, for every geodesically convex function $\varphi:\mathcal M\to\mathbb R\cup\{+\infty\}$,
			\begin{align*}
				|\nabla\varphi|(\bar u)<\gamma\quad\implies\quad \argmin_{v\in\mathcal M}\big\{\mathcal J(v) + \varphi(v)\big\}\cap\mathbb B(\bar u,\delta)\subseteq \{\bar u\}.
			\end{align*}
		\end{itemize}
	\end{corollary}
\noindent 
	An analogous result as the one of Proposition \ref{P1} can also be obtained for geodesically convex functionals when the local slope is considered. 
	\begin{proposition}\label{P2}
		Suppose that $(\mathcal M, d)$ is a complete geodesic metric space and $\mathcal J:\mathcal M\to\mathbb R\cup\{+\infty\}$ is geodesically convex. Let $\gamma>0$. The following statements are equivalent. 
        \begin{itemize}
            \item[$(i)$] There exists $\delta\in(0,+\infty]$ such that 
            \begin{align*}
			|\nabla \mathcal J|(u)\ge\gamma\quad\forall  u\in \mathbb B(\bar u,\delta)\setminus\{\bar u\}.
		\end{align*}

            \item[$(ii)$] There exists $\delta\in(0,+\infty]$ such that $\bar u$ is a sharp local    minimizer with parameters $\delta$ and $\gamma$.
        \end{itemize}
	\end{proposition}

	\section{Proofs}\label{proofs}

	\subsection{Proof of Theorem \ref{Thm1}} 
	We begin with implication {$(iii)\implies(i)$}. 
	Let $\xi\in X^*$ with $\|\xi\|<\gamma$. Then, 
	\begin{align*}
		f(x) - \langle \xi, x\rangle &\ge f(\bar x) +\gamma \|x-\bar x\| -  \langle \xi, x\rangle  \ge f(\bar x) -\langle \xi, \bar x \rangle + (\gamma - \|\xi\|)\|x-\bar x\|.
	\end{align*}
	Thus, $\bar x$ is a strict minimizer of $f - \langle \xi, \cdot\rangle$, whence the result follows. 
	The implication {$(i)\implies(ii)$} is trivial. Let us then proceed with implication  {$(ii)\implies(iii)$}.
	Let $\Theta\subseteq X^*$ dense and  $\gamma>0$ be such that
	\begin{align}\label{dens}
		\bar x\in \argmin_{y\in X}\big\{f(y) -\langle \xi, y\rangle\big\}\quad \text{for all $\xi\in \Theta$ satisfying $\|\xi\|<\gamma$.}
	\end{align}
	Let $x\in X$ and $\varepsilon>0$ be given. By the hyperplane separation theorem (geometric Hahn-Banach), we can find  $\theta\in X^*\setminus\{0\}$ such that 
	\begin{align}\label{ghb}
		\sup_{y\in\mathbb B(x,\|x-\bar x\|)} \langle \theta,y \rangle\le \langle \theta, \bar x\rangle 
	\end{align}
	Using the definition of dual norm, a straightforward calculation yields
	\begin{equation}\label{supdn}
		\begin{aligned}
			\sup_{y\in\mathbb B(x,\|x-\bar x\|)}\langle \theta,y \rangle = \sup_{y\in\mathbb B(x,\|x-\bar x\|)}\langle \theta,y - x \rangle + \langle \theta ,x \rangle = \|\theta\|\|x-\bar x\| +  \langle\theta,x\rangle.
		\end{aligned}
	\end{equation}
	Define $\xi:=-(\gamma-\varepsilon)\theta/\|\theta\|$. Observe that, by (\ref{ghb}) and  (\ref{supdn}), 
	\begin{equation}\label{ghb2}
		\begin{aligned}
			\langle \xi,x-\bar x\rangle = 	\frac{\gamma-\varepsilon}{\|\theta\|}\langle \theta,\bar x\rangle -\frac{\gamma-\varepsilon}{\|\theta\|}	\langle \theta,x\rangle \ge \frac{\gamma-\varepsilon}{\|\theta\|} \|\theta\|\|x-\bar x\|  = (\gamma-\varepsilon)\|x-\bar x\|.
		\end{aligned}
	\end{equation}
	Let $\{\xi_n\}_{n\in\mathbb N}\subseteq\Theta$ be such that $ \|\xi_n - \xi\| \longrightarrow0.$
	Since $\|\xi\|<\gamma$, there exists $N\in \mathbb N$ such that  $\|\xi_n\|<\gamma$ for all $n\ge N$. Then, by assumption (\ref{dens}), there holds $\bar x\in \argmin_{y\in X}\big\{f(y)-\langle \xi_n,y\rangle\}$ for all $n\ge N$. This can be written as 
	\begin{align*}
		f(\bar x) -\langle \xi_n ,\bar x\rangle \le  f(x) -\langle \xi_n , x\rangle \quad \forall n\ge N. 
	\end{align*}
	Taking limit as $n\to+\infty$, we get $f(\bar x) -\langle \xi ,\bar x\rangle \le  f(x) -\langle \xi , x\rangle$. From this, 
	\begin{align*}
		f(x)-f(\bar x)\ge \langle \xi, x-\bar x\rangle = (\gamma - \varepsilon) \|x-\bar x\|.
	\end{align*}
	Since $x\in X$ and $\varepsilon>0$ were arbitrary, we conclude that 
	\begin{align*}
		\inf_{x\in X\setminus\{\bar x\}}\frac{f(x)-f(\bar x)}{\|x-\bar x\|}\ge\gamma -\varepsilon\quad \forall \varepsilon>0. 
	\end{align*}
	Letting $\varepsilon\longrightarrow0^+$ yields the result.  \hfill\(\square\)

	\subsection{Proof of Corollary \ref{Cor1}}
	Implication $(ii)\implies (i)$ is trivial because bounded linear functionals $\xi:X\to\mathbb R$ are Lipschitz. Let  $\zeta:X\to\mathbb R$ be a Lipschitz function with $\Lip \zeta<\gamma$.
	By Theorem \ref{Thm1}, $\bar x$ is a sharp minimizer of $f$ and consequently
	\begin{align*}
		f(x) +\zeta (x) &\ge f(\bar x) +\gamma \|x-\bar x\| +\zeta(x) \ge f(\bar x) + \zeta(\bar x) + (\gamma - \Lip\zeta)\|x-\bar x\|
	\end{align*}
	for all $x\in X$. We conclude that $\bar x$ is the unique minimizer of $f+\zeta$. \hfill\(\square\)

	\subsection{Proof of Proposition \ref{P1}}\label{pss}
	Suppose that $\bar x\in X$ is a sharp minimizer with modulus $\gamma>0$. Suppose there exists $x\in X\setminus\{\bar x\}$ with $|\partial f|(x)<\gamma$, then 
	\begin{align*}
		f(x) - f(y) < \gamma \|x-y\|\quad \forall y\in X\setminus\{x\}.
	\end{align*}
	Combining this with the fact that $\bar x$ is a sharp minimizer, we get
	\begin{align*}
		\gamma\|x-\bar x\| \le f(x) - f(\bar x)< \gamma\|x-\bar x\|.
	\end{align*}
	This is a contradiction, and hence $|\partial f|(x)\ge\gamma$. 
	\smallbreak\noindent  
	Suppose now that there holds
	\begin{align}\label{assinle}
		|\partial f|(y) \ge\gamma \quad \forall y\in X\setminus\{\bar x\}. 
	\end{align}
    We first claim that $\bar x$ is a global minimizer of $f$. 
Fix $\eta\in(0,\gamma)$ and $\varepsilon>0$. Choose $x_0\in X$ such that
$f(x_0)\le \inf_{y\in X} f(y)+\varepsilon$.
By Ekeland's variational principle applied to $f$ with parameters
$\varepsilon$ and $\lambda:=\varepsilon/\eta$, there exists $x_\varepsilon\in X$ such that
$f(x_\varepsilon)\le f(x_0)$ and
\begin{align}\label{ek_min}
f(x_\varepsilon)\le f(y)+\eta\|y-x_\varepsilon\|\qquad \forall y\in X.
\end{align}
Let $y\in X\setminus\{x_\varepsilon\}$. From \eqref{ek_min} we obtain
\[
\frac{\max\{f(x_\varepsilon)-f(y),0\}}{\|x_\varepsilon-y\|}
\le \eta,
\]
hence $|\partial f|(x_\varepsilon)\le \eta<\gamma$. By \eqref{assinle} it follows that
$x_\varepsilon=\bar x$. Therefore
\[
f(\bar x)=f(x_\varepsilon)\le f(x_0)\le \inf_{y\in X} f(y)+\varepsilon.
\]
Since $\varepsilon>0$ is arbitrary, we conclude $f(\bar x)=\inf_{y\in X} f(y)$, proving the claim.
\smallbreak\noindent
	Now, let $\gamma'\in(0,\gamma)$ be arbitrary. Suppose there exists $x\in X\setminus\{\bar x\}$ such that
	\begin{align}\label{inco}
		f(x)<f(\bar x) + {\gamma}'  \|x-\bar x\|.
	\end{align}
	By Ekeland's variational principle, there exists $\hat x\in X$ such that
	\begin{align*}
		\| x-\hat x\| \le \frac{2\gamma'}{\gamma'+\gamma}\hspace*{0.02cm}\|x-\bar x\| \quad \text{and}\quad  \hat x\in \argmin_{y\in X}\big\{f(y) + \frac{\gamma+\gamma'}{2}\| y-\hat x\|\}.
	\end{align*}
	We can conclude then that
	\begin{align*}
		\frac{f(\hat x)-f(y)}{\|\hat x-y\|} \le \frac{\gamma+\gamma'}{2}\quad \forall y\in X\setminus\{\hat x\}.
	\end{align*}
	Hence, $|\partial f|(\hat x)\le 2^{-1}(\gamma+\gamma')<\gamma$; so by (\ref{assinle}), $\hat x=\bar x$. From this, $\|x-\bar x\|\le 2\gamma'(\gamma+\gamma')^{-1} \|x-\bar x\|<\|x-\bar x\|$; a contradiction to (\ref{inco}). We conclude that 
	\begin{align*}
		f(x) \ge f(\bar x) + {\gamma}' \|x-\bar x\| \quad \forall x\in X.
	\end{align*}
	The result follows since $\gamma'\in(0,\gamma)$ was arbitrary. \hfill\(\square\)

    \subsection{Proof of Proposition \ref{pappendix}}
    It is easy to check that $f^{**}(\bar x)=f(\bar x)$ and $\bar x\in \argmin_{y\in X} f^{**}(y)$. 
Assume now that $\bar x$ is a sharp minimizer of $f$ with modulus $\gamma$. By Theorem~\ref{Thm1}, for every $\xi\in X^*$ with $\|\xi\|<\gamma$, the point $\bar x$ minimizes $f-\langle \xi,\cdot\rangle$. Equivalently,
$f^*(\xi)=\langle \xi,\bar x\rangle-f(\bar x)$ whenever $\|\xi\|<\gamma$.
Thus, for every $x\in X$,
\begin{align*}
f^{**}(x)
&\ge \sup_{\|\xi\|<\gamma}\bigl\{\langle \xi,x\rangle-f^*(\xi)\bigr\}=f(\bar x)+\sup_{\|\xi\|<\gamma}\langle \xi,x-\bar x\rangle=f(\bar x)+\gamma\|x-\bar x\|.
\end{align*}
Using that $f^{**}(\bar x)=f(\bar x)$ this becomes $f^{**}(x)\ge f^{**}(\bar x)+\gamma\|x-\bar x\|$ for all $x\in X$. We conclude that $\bar x$ is a sharp minimizer of $f^{**}$ with modulus $\gamma$. 
Conversely, assume that $\bar x$ is a sharp minimizer of $f^{**}$ with modulus $\gamma$. Then
$f^{**}(x)\ge f^{**}(\bar x)+\gamma\|x-\bar x\|$ for every $x\in X$.
Since $f^{**}\le f$ and $f^{**}(\bar x)=f(\bar x)$, we get
\[
f(x)\ge f^{**}(x)\ge f(\bar x)+\gamma\|x-\bar x\|\quad \forall x\in X.
\]
Hence $\bar x$ is a sharp minimizer of $f$ with modulus $\gamma$. \hfill\(\square\)

	\subsection{Proof of Theorem \ref{Thm2}}
	Suppose $\bar u\in\mathcal M$ is a sharp local    minimizer with parameters $\delta$ and $\gamma$. Let $\zeta:\mathcal M\to\mathbb R$ be a Lipschitz function with $\Lip\zeta<\gamma$. 
	Let $u\in\argmin_{v\in\mathcal M}\{\mathcal J(v)+\zeta(v)\}\cap\mathbb B(\bar u,\delta)$; then $\mathcal J(u)+\zeta(u)\le \mathcal J(\bar u)+\zeta(\bar u)$. 
	Since $\bar u$ is a sharp local    minimizer, 
	\begin{align*}
		\gamma d(u,\bar u)\le \mathcal J(u)-\mathcal J(\bar u)\le \zeta(\bar u) - \zeta(u) \le \Lip \zeta\, d(u,\bar u).
	\end{align*}
	If \( u\neq \bar u \), then \( \gamma\le \Lip\zeta \); a contradiction. We have proved that 
	\begin{align}\label{genincl}
		\argmin_{v\in \mathcal M}\{\mathcal{J}(v)+\zeta(v)\}\cap \mathbb B(\bar u,\delta)\subseteq \{\bar u\}.
	\end{align}
	Assume now that there exists $\delta\in(0,+\infty]$ such that  (\ref{genincl}) holds for any Lipschitz function $\zeta:\mathcal M\to \mathbb R$. Let $\gamma'<\gamma$ be given.  Suppose that there exists $u\in\mathbb B(\bar u,\delta/2)\setminus\{\bar u\}$ such that 
	\begin{align*}
		\mathcal J(u) < \mathcal J(\bar u) + {\gamma}' d(u,\bar u).
	\end{align*} 
	By Ekeland's variational principle, there exists $\hat u\in \mathcal M$ such that 
	\begin{itemize}
		\item[$(a)$] $\mathcal J(\hat u)\le \mathcal J(u)$;
		
		\item[$(b)$] $\displaystyle d(u,\hat u)\le \frac{2\gamma'}{\gamma'+\gamma} d(u,\bar u)$;
		
		\item[$(c)$] $\displaystyle \argmin_{v\in\mathcal M}\big\{\mathcal J(v) + \frac{\gamma'+\gamma}{2} d(v,\hat u)  \big\}=\{\hat u\}$
	\end{itemize}
	Let $\zeta:\mathcal M\to\mathbb R$ be given by $\zeta(v):=2^{-1}(\gamma'+\gamma)d(v,\hat u)$. Observe that
    \[
        d(\hat u,\bar u)\le d(\hat u,u) + d(u,\bar u)\le \Big(\frac{2\gamma'}{\gamma'+\gamma}  +1\Big)d(u,\bar u) \le 2d(u,\bar u)\le \delta.
    \]
    We see then that $\hat u\in\argmin_{v\in \mathcal M}\{\mathcal{J}(v)+\zeta(v)\}\cap \mathbb B(\bar u,\delta)$. Combining this with the fact that   $\Lip\zeta=2^{-1}(\gamma'+\gamma)<\gamma$, we get, by hypothesis, that $\hat u=\bar u$. This is a contradiction to item $(b)$ since $2\gamma'(\gamma'+\gamma)^{-1}<1$. We conclude that 
	\begin{align*}
		\mathcal J(u) \ge \mathcal J(\bar u) + {\gamma}'d(u,\bar u)\quad \forall u\in\mathbb B(\bar u,\delta/2).
	\end{align*} 
	Since $\gamma'<\gamma$ was arbitrary, the result follows. \hfill\(\square\)
	
	\subsection{Proof of Corollary \ref{Cor2}}
	Let us begin proving implication $(i)\implies(ii)$. 
	 Let $\varphi:\mathcal M\to\mathbb R\cup\{+\infty\}$ be a geodesically convex functional satisfying $|\nabla \varphi|(\bar u)<\gamma$. From Theorem \ref{Thm2}, $\bar u$ must be a sharp local    minimizer with parameters $\delta'$ and $\gamma$ for some $\delta'>0$. Then, using \cite[Theorem 2.4.9]{AGS_2008}, for any $u\in\argmin_{v\in\mathcal M}\{\mathcal J(v)+\varphi(v)\}\cap\mathbb B(\bar u,\delta')$,
	\begin{align*}
		\gamma d(u,\bar u)\le\mathcal J(u)-\mathcal J(\bar u) \le \varphi(\bar u)-\varphi(u) \le |\nabla\varphi|(\bar u) d(u,\bar u).
	\end{align*}
	We conclude that if $u\neq \bar u$, then $|\nabla\varphi|(\bar u)\ge \gamma$; a contradiction, hence $u=\bar u$. 
	\smallbreak 
	We prove now implication $(ii)\implies(i)$. 
	Let $\zeta:\mathcal M\to\mathbb R$ be an arbitrary Lipschitz function with $\Lip\zeta<\gamma$.  Fix
	\begin{align*}
		u\in\argmin_{v\in\mathcal M}\{\mathcal J(v)+\zeta(v)\}\cap\mathbb B(\bar u,\delta).
	\end{align*}
	Then, $	\mathcal J(u) + \zeta(u) \le \mathcal J(v) + \zeta(v)$ for all $v\in\mathcal M$. From this, 
	\begin{align}\label{someineq}
		\mathcal J(u) -\mathcal J(v) \le \zeta(v) - \zeta(u) \le \Lip\zeta\, d(v,u)\quad \forall v\in\mathcal M. 
	\end{align}
	Define $\varphi:\mathcal M\to\mathbb R$ by $\varphi(v):=\Lip\zeta\,d(v,u)$. Then, by (\ref{someineq}), 
	\begin{align*}
		\mathcal J(u) + \varphi(u) = \mathcal J(u)\le \mathcal J(v) + \Lip \zeta\, d(u,v)=\mathcal J(v)+\varphi(v)\quad \forall v\in\mathcal M. 
	\end{align*}
	We see that $u\in\argmin_{v\in\mathcal M}\{\mathcal J(v)+\varphi(v)\}\cap\mathbb B(\bar u,\delta)$. Since $\varphi$ is geodesically convex and $|\nabla\varphi|(\bar u)\le\Lip\zeta<\gamma$, we conclude $u=\bar u$. \hfill\(\square\)

	\subsection{Proof of Proposition \ref{P2}}
	Assume 	$|\nabla \mathcal J|(u)\ge\gamma$ for all $u\in \mathbb B(\bar u,\delta)\setminus\{\bar u\}.$ Let $\zeta:\mathcal M\to\mathbb R$ be a Lipschitz function satisfying $\Lip\zeta<\gamma$. We will prove that $\argmin_{v\in\mathcal M}\{\mathcal J(v)+\zeta(v)\}\cap\mathbb B(\bar u,\delta)\subseteq\{\bar u\}$, and then by Theorem \ref{Thm2}, we will conclude that $\bar u$ is a sharp local    minimizer. 
	Let $u\in\argmin_{v\in \mathcal M}\{\mathcal J(v)+\zeta(v)\}\cap\mathbb B(\bar u,\delta)$, then $\mathcal J(u)+\zeta(u)\le \mathcal J(v)+\zeta(v)$ for all $v\in\mathcal M$; this  implies $\mathcal J(u)-\mathcal J(v)\le \Lip\zeta\,d(u,v)$ for all $v\in\mathcal M$. Hence, 
	\begin{align*}
		\max\big\{\mathcal J(u)-\mathcal J(v),0\big\}\le \Lip\zeta\, d(u,v)\quad \forall v\in\mathcal M.  
	\end{align*}
	From this follows that $|\nabla\mathcal J|(u)\le \Lip \zeta<\gamma$. Hence, by assumption $u=\bar u$. 
	\smallbreak\noindent 
	Now assume that $\bar u$ is a sharp local    minimizer with parameters $\delta\in(0,+\infty]$ and $\gamma>0$. Let $u\in\mathbb B(\bar u,\delta)\setminus\{\bar u\}$. Since $\mathcal J$ is geodesically convex,  $\mathcal J(u)-\mathcal J(v)\le |\nabla\mathcal J|(u)\,d(u,v)$ for all $v\in \mathcal M$, see \cite[Theorem 2.4.9]{AGS_2008}. In particular, 
	\begin{align*}
		\gamma\hspace*{0.04cm} d(u,\bar u)\le \mathcal J(u) - \mathcal J(\bar u) \le |\nabla\mathcal J|(u)\,d(u,\bar u).
	\end{align*}
	Whence we  conclude $\gamma\le|\nabla\mathcal J|(u)$. \hfill\(\square\)

\appendix
\section{Nonpositively curved spaces and geodesic convexity}\label{App:CAT0}
\subsection*{Geodesics and convexity}\label{App:geodesics}
Throughout this subsection, we consider a metric space $(\mathcal M,d)$.  
A continuous curve $\sigma:[0,1]\to\mathcal M$ is called a (constant-speed) geodesic if
\begin{equation}\label{eq:geodesic_def}
d\big(\sigma(t),\sigma(s)\big)=|t-s|\,d\big(\sigma(0),\sigma(1)\big)\qquad \forall s,t\in[0,1].
\end{equation}
We say that two points $u,v\in\mathcal M$ can be joined by a geodesic if there exists a geodesic $\sigma:[0,1]\to\mathcal M$ such that $\sigma(0)=u$ and $\sigma(1)=v$; in particular, one has
\begin{equation}\label{eq:geodesic_endpoints}
d\big(\sigma(t),u\big)=t\,d(u,v)\quad \text{and\quad}d\big(\sigma(t),v\big)=(1-t)\,d(u,v)\qquad \forall t\in[0,1].
\end{equation}
We say that $(\mathcal M,d)$ is a geodesic metric space if for every pair of points $u,v\in\mathcal M$ there exists at least one geodesic joining $u$ and $v$.
\smallbreak\noindent
If $\sigma:[0,1]\to\mathcal M$ is a geodesic, then its metric derivative exists at every $t\in(0,1)$ and is given by
\[
|\dot\sigma|(t):=\lim_{h\to 0}\frac{d\big(\sigma(t+h),\sigma(t)\big)}{|h|}
= d\big(\sigma(0),\sigma(1)\big).
\]
This explains the use of the term \emph{constant-speed} for such geodesics in metric spaces.
\smallbreak\noindent
We say that a proper functional $\mathcal J:\mathcal M\to\mathbb R\cup\{+\infty\}$ is \emph{geodesically convex} if for every $u,v\in \mathcal M$ there exists a geodesic $\sigma:[0,1]\to\mathcal M$ joining $u$ and $v$ such that
\begin{equation}\label{eq:gconv}
\mathcal J(\sigma(t))\le (1-t)\,\mathcal J(u)+t\,\mathcal J(v)\qquad \forall t\in[0,1].
\end{equation}
\subsection*{Nonpositive curvature}\label{App:cat0}

A geodesic metric space $(\mathcal M,d)$ is called a $\operatorname{CAT}(0)$ space if for every $u,v,w\in\mathcal M$ there exists a geodesic $\sigma:[0,1]\to\mathcal M$ with $\sigma(0)=v$ and $\sigma(1)=w$ such that
\begin{equation}\label{eq:cat0}
d\big(\sigma(t),u\big)^2
\le (1-t)\,d(v,u)^2+t\,d(w,u)^2-t(1-t)\,d(v,w)^2
\qquad \forall t\in[0,1].
\end{equation}
A complete $\operatorname{CAT}(0)$ space is often referred to as an \emph{Hadamard space}.
 It is easy to see directly from \eqref{eq:cat0}, that in a $\operatorname{CAT}(0)$ space any two points can be joined by a unique geodesic.
\smallbreak\noindent
Inequality \eqref{eq:cat0} can be viewed as a non-flat metric analogue of the Euclidean identity for straight line segments. It can be read as follows: if one moves from one point to another along a geodesic, then the squared distance to an arbitrary base point cannot exceed the value predicted by the Euclidean model. In particular, geodesics do not diverge faster than in the Euclidean comparison geometry, reflecting the “thin triangle” behavior characteristic of nonpositive curvature.
\smallbreak\noindent
A relevant feature of nonpositively curved spaces is that distance functions are geodesically convex; i.e.,  if $(\mathcal M,d)$ is $\operatorname{CAT}(0)$, then for every fixed $u\in\mathcal M$ the map $v\mapsto d(v,u)$ is geodesically convex.

\section*{Acknowledgments}
The author expresses gratitude to Tr\'i Minh L\^e for fruitful discussions on sharp minimizers. 

\bibliography{references}{}

@book {P_1987,
	AUTHOR = {Polyak, Boris T.},
	TITLE = {Introduction to optimization},
	SERIES = {Translations Series in Mathematics and Engineering},
	NOTE = {Translated from the Russian,
	With a foreword by Dimitri P. Bertsekas},
	PUBLISHER = {Optimization Software, Inc., Publications Division, New York},
	YEAR = {1987},
	PAGES = {xxvii+438},
	ISBN = {0-911575-14-6},
	MRCLASS = {49-01 (65Kxx 90Cxx)},
	MRNUMBER = {1099605},
}

@phdthesis{F_1988,
	title={Weak sharp minima and penalty functions in mathematical programming.},
	author={Ferris, Michael Charles},
	year={1988},
	school={University of Cambridge}
}

@article {BF_1993,
	AUTHOR = {Burke, J. V. and Ferris, M. C.},
	TITLE = {Weak sharp minima in mathematical programming},
	JOURNAL = {SIAM J. Control Optim.},
	FJOURNAL = {SIAM Journal on Control and Optimization},
	VOLUME = {31},
	YEAR = {1993},
	NUMBER = {5},
	PAGES = {1340--1359},
	ISSN = {0363-0129},
	MRCLASS = {90C31 (65K05)},
	MRNUMBER = {1234006},
	MRREVIEWER = {Doug\ Ward},
	DOI = {10.1137/0331063},
	URL = {https://doi.org/10.1137/0331063},
}

@article {U_2015,
	AUTHOR = {Uderzo, A.},
	TITLE = {On the variational behaviour of functions with positive
	steepest descent rate},
	JOURNAL = {Positivity},
	FJOURNAL = {Positivity. An International Mathematics Journal Devoted to
	Theory and Applications of Positivity},
	VOLUME = {19},
	YEAR = {2015},
	NUMBER = {4},
	PAGES = {725--745},
	ISSN = {1385-1292,1572-9281},
	MRCLASS = {49K40 (49J52 90C25 90C31)},
	MRNUMBER = {3415100},
	MRREVIEWER = {Tran\ T. A. Nghia},
	DOI = {10.1007/s11117-015-0324-x},
	URL = {https://doi.org/10.1007/s11117-015-0324-x},
}

@book {AGS_2008,
	AUTHOR = {Ambrosio, Luigi and Gigli, Nicola and Savar\'e, Giuseppe},
	TITLE = {Gradient flows in metric spaces and in the space of
	probability measures},
	SERIES = {Lectures in Mathematics ETH Z\"urich},
	EDITION = {Second},
	PUBLISHER = {Birkh\"auser Verlag, Basel},
	YEAR = {2008},
	PAGES = {x+334},
	ISBN = {978-3-7643-8721-1},
	MRCLASS = {49-02 (28A33 35K55 35K90 49Q20 60B05)},
	MRNUMBER = {2401600},
	MRREVIEWER = {Pietro\ Celada},
}

@article {L_2009,
	AUTHOR = {Lassonde, Marc},
	TITLE = {Asplund spaces, {S}tegall variational principle and the {RNP}},
	JOURNAL = {Set-Valued Var. Anal.},
	FJOURNAL = {Set-Valued and Variational Analysis. Theory and Applications},
	VOLUME = {17},
	YEAR = {2009},
	NUMBER = {2},
	PAGES = {183--193},
	ISSN = {1877-0533,1877-0541},
	MRCLASS = {46B22 (46B10 46B26 46N10 49J50)},
	MRNUMBER = {2529695},
	MRREVIEWER = {Mari\'an\ Fabian},
	DOI = {10.1007/s11228-009-0111-6},
	URL = {https://doi.org/10.1007/s11228-009-0111-6},
}

@inproceedings{KMNST_2017,
	author = {Keskar, Nitish Shirish and Mudigere, Dheevatsa and Nocedal, Jorge and Smelyanskiy, Misha and Tang, Ping Tak Peter},
	title = {On Large-Batch Training for Deep Learning: Generalization Gap and Sharp Minima},
	booktitle = {Proceedings of the International Conference on Learning Representations (ICLR)},
	year = {2017}
}

@inproceedings{DPBB_2017,
	title={Sharp minima can generalize for deep nets},
	author={Dinh, Laurent and Pascanu, Razvan and Bengio, Samy and Bengio, Yoshua},
	booktitle={International Conference on Machine Learning},
	pages={1019--1028},
	year={2017}
}

@article {BF_1995,
	AUTHOR = {Burke, J. V. and Ferris, M. C.},
	TITLE = {A {G}auss-{N}ewton method for convex composite optimization},
	JOURNAL = {Math. Programming},
	FJOURNAL = {Mathematical Programming},
	VOLUME = {71},
	YEAR = {1995},
	NUMBER = {2},
	PAGES = {179--194},
	ISSN = {0025-5610,1436-4646},
	MRCLASS = {90C30 (65K05)},
	MRNUMBER = {1373362},
	MRREVIEWER = {Asen\ L.\ Dontchev},
	DOI = {10.1007/BF01585997},
	URL = {https://doi.org/10.1007/BF01585997},
}

@article {FNR_2023,
	AUTHOR = {Fadili, Jalal and Nghia, Tran T. A. and Tran, Trinh T. T.},
	TITLE = {Sharp, strong and unique minimizers for low complexity robust
	recovery},
	JOURNAL = {Inf. Inference},
	FJOURNAL = {Information and Inference. A Journal of the IMA},
	VOLUME = {12},
	YEAR = {2023},
	NUMBER = {3},
	PAGES = {Paper No. iaad005, 53},
	ISSN = {2049-8764,2049-8772},
	MRCLASS = {90C26 (49N45)},
	MRNUMBER = {4581451},
	MRREVIEWER = {Ji\ Li},
	DOI = {10.1093/imaiai/iaad005},
	URL = {https://doi.org/10.1093/imaiai/iaad005},
}

@article{BD_2002,
	AUTHOR = {Burke, James V. and Deng, Sien},
	TITLE = {Weak sharp minima revisited. {I}. {B}asic theory},
	JOURNAL = {Control Cybernet.},
	FJOURNAL = {Control and Cybernetics},
	VOLUME = {31},
	YEAR = {2002},
	NUMBER = {3},
	PAGES = {439--469},
	ISSN = {0324-8569,2720-4278},
	MRCLASS = {90C31 (49J52)},
	MRNUMBER = {1978735},
}

@article {BD_2005,
	AUTHOR = {Burke, James V. and Deng, Sien},
	TITLE = {Weak sharp minima revisited. {II}. {A}pplication to linear
	regularity and error bounds},
	JOURNAL = {Math. Program.},
	FJOURNAL = {Mathematical Programming. A Publication of the Mathematical
	Programming Society},
	VOLUME = {104},
	YEAR = {2005},
	NUMBER = {2-3},
	PAGES = {235--261},
	ISSN = {0025-5610,1436-4646},
	MRCLASS = {90C31 (49J52)},
	MRNUMBER = {2179237},
	DOI = {10.1007/s10107-005-0615-2},
	URL = {https://doi.org/10.1007/s10107-005-0615-2},
}

@article {BD_2009,
	AUTHOR = {Burke, James V. and Deng, Sien},
	TITLE = {Weak sharp minima revisited. {III}. {E}rror bounds for
	differentiable convex inclusions},
	JOURNAL = {Math. Program.},
	FJOURNAL = {Mathematical Programming. A Publication of the Mathematical
	Programming Society},
	VOLUME = {116},
	YEAR = {2009},
	NUMBER = {1-2},
	PAGES = {37--56},
	ISSN = {0025-5610,1436-4646},
	MRCLASS = {90C31 (49J52)},
	MRNUMBER = {2421272},
	DOI = {10.1007/s10107-007-0130-8},
	URL = {https://doi.org/10.1007/s10107-007-0130-8},
}

@article {KZ_2003,
	AUTHOR = {Ng, Kung Fu and Zheng, Xi Yin},
	TITLE = {Global weak sharp minima on {B}anach spaces},
	JOURNAL = {SIAM J. Control Optim.},
	FJOURNAL = {SIAM Journal on Control and Optimization},
	VOLUME = {41},
	YEAR = {2003},
	NUMBER = {6},
	PAGES = {1868--1885},
	ISSN = {0363-0129,1095-7138},
	MRCLASS = {49J52 (90C31)},
	MRNUMBER = {1972538},
	MRREVIEWER = {Doug\ Ward},
	DOI = {10.1137/S0363012901389469},
	URL = {https://doi.org/10.1137/S0363012901389469},
}

@article {U_2014,
	AUTHOR = {Uderzo, A.},
	TITLE = {Some dual conditions for global weak sharp minimality of
	nonconvex functions},
	JOURNAL = {Numer. Funct. Anal. Optim.},
	FJOURNAL = {Numerical Functional Analysis and Optimization. An
	International Journal},
	VOLUME = {35},
	YEAR = {2014},
	NUMBER = {7-9},
	PAGES = {1258--1285},
	ISSN = {0163-0563,1532-2467},
	MRCLASS = {49J53 (49K40 90C26 90C48)},
	MRNUMBER = {3230101},
	MRREVIEWER = {Doug\ Ward},
	DOI = {10.1080/01630563.2014.895749},
	URL = {https://doi.org/10.1080/01630563.2014.895749},
}

@article {ZX_2011,
	AUTHOR = {Zhou, Jinchuan and Xu, Xiuhua},
	TITLE = {Equivalent properties of global weak sharp minima with
	applications},
	JOURNAL = {J. Inequal. Appl.},
	FJOURNAL = {Journal of Inequalities and Applications},
	YEAR = {2011},
	PAGES = {2011:137, 9},
	ISSN = {1029-242X},
	MRCLASS = {90C30 (90C25 90C26 90C34)},
	MRNUMBER = {2887135},
	MRREVIEWER = {Jin\ Ling\ Zhao},
	DOI = {10.1186/1029-242X-2011-137},
	URL = {https://doi.org/10.1186/1029-242X-2011-137},
}

@article {W_1994,
	AUTHOR = {Ward, D. E.},
	TITLE = {Characterizations of strict local minima and necessary
	conditions for weak sharp minima},
	JOURNAL = {J. Optim. Theory Appl.},
	FJOURNAL = {Journal of Optimization Theory and Applications},
	VOLUME = {80},
	YEAR = {1994},
	NUMBER = {3},
	PAGES = {551--571},
	ISSN = {0022-3239,1573-2878},
	MRCLASS = {90C31 (49J52 49K40)},
	MRNUMBER = {1265176},
	MRREVIEWER = {Nguyen Van Thoai},
	DOI = {10.1007/BF02207780},
	URL = {https://doi.org/10.1007/BF02207780},
}

@article {SW_1999,
	AUTHOR = {Studniarski, Marcin and Ward, Doug E.},
	TITLE = {Weak sharp minima: characterizations and sufficient
	conditions},
	JOURNAL = {SIAM J. Control Optim.},
	FJOURNAL = {SIAM Journal on Control and Optimization},
	VOLUME = {38},
	YEAR = {1999},
	NUMBER = {1},
	PAGES = {219--236},
	ISSN = {0363-0129,1095-7138},
	MRCLASS = {90C31 (49J52)},
	MRNUMBER = {1740599},
	MRREVIEWER = {Alexander\ Shapiro},
	DOI = {10.1137/S0363012996301269},
	URL = {https://doi.org/10.1137/S0363012996301269},
}

@article {ZMX_2012,
	AUTHOR = {Zhou, Jinchuan and Mordukhovich, Boris S. and Xiu, Naihua},
	TITLE = {Complete characterizations of local weak sharp minima with
	applications to semi-infinite optimization and
	complementarity},
	JOURNAL = {Nonlinear Anal.},
	FJOURNAL = {Nonlinear Analysis. Theory, Methods \& Applications. An
	International Multidisciplinary Journal},
	VOLUME = {75},
	YEAR = {2012},
	NUMBER = {3},
	PAGES = {1700--1718},
	ISSN = {0362-546X,1873-5215},
	MRCLASS = {49J52 (49K40 90C26)},
	MRNUMBER = {2861368},
	DOI = {10.1016/j.na.2011.05.084},
	URL = {https://doi.org/10.1016/j.na.2011.05.084},
}

@article {ZW_2012,
	AUTHOR = {Zhou, Jinchuan and Wang, Changyu},
	TITLE = {New characterizations of weak sharp minima},
	JOURNAL = {Optim. Lett.},
	FJOURNAL = {Optimization Letters},
	VOLUME = {6},
	YEAR = {2012},
	NUMBER = {8},
	PAGES = {1773--1785},
	ISSN = {1862-4472,1862-4480},
	MRCLASS = {90C31 (49K40)},
	MRNUMBER = {2996483},
	MRREVIEWER = {Hongxia\ Yin},
	DOI = {10.1007/s11590-011-0369-0},
	URL = {https://doi.org/10.1007/s11590-011-0369-0},
}

@article {LT_2024,
	AUTHOR = {L\^e, Tr\'i Minh and Tapia-Garc\'ia, Sebasti\'an},
	TITLE = {On (discounted) global {E}ikonal equations in metric spaces},
	JOURNAL = {arXiv preprint arXiv:2410.00530},
	YEAR = {2024},
	URL = {https://arxiv.org/abs/2410.00530}
}

@article {LSZ_2021,
    AUTHOR = {Liu, Qing and Shanmugalingam, Nageswari and Zhou, Xiaodan},
     TITLE = {Equivalence of solutions of eikonal equation in metric spaces},
   JOURNAL = {J. Differential Equations},
  FJOURNAL = {Journal of Differential Equations},
    VOLUME = {272},
      YEAR = {2021},
     PAGES = {979--1014},
      ISSN = {0022-0396,1090-2732},
   MRCLASS = {35F20 (30F30 35D40 35F21 49L25)},
  MRNUMBER = {4166072},
       DOI = {10.1016/j.jde.2020.10.018},
       URL = {https://doi.org/10.1016/j.jde.2020.10.018},
}

@article {GHN_2015,
	AUTHOR = {Giga, Yoshikazu and Hamamuki, Noboru and Nakayasu, Akitoshi},
	TITLE = {Eikonal equations in metric spaces},
	JOURNAL = {Trans. Amer. Math. Soc.},
	FJOURNAL = {Transactions of the American Mathematical Society},
	VOLUME = {367},
	YEAR = {2015},
	PAGES = {49--66},
	ISSN = {0002-9947},
	DOI = {10.1090/S0002-9947-2014-06048-4},
	URL = {https://doi.org/10.1090/S0002-9947-2014-06048-4}
}

@article {GS_2015,
	AUTHOR = {Gangbo, Wilfrid and \'Swi{e}ch, Andrzej},
	TITLE = {Metric viscosity solutions of {H}amilton--{J}acobi equations depending on local slopes},
	JOURNAL = {Calc. Var. Partial Differential Equations},
	FJOURNAL = {Calculus of Variations and Partial Differential Equations},
	VOLUME = {54},
	YEAR = {2015},
	PAGES = {1183--1218},
	ISSN = {0944-2669},
	DOI = {10.1007/s00526-015-0852-2},
	URL = {https://doi.org/10.1007/s00526-015-0852-2}
}

@article {DS_2022,
    AUTHOR = {Daniilidis, Aris and Salas, David},
     TITLE = {A determination theorem in terms of the metric slope},
   JOURNAL = {Proc. Amer. Math. Soc.},
  FJOURNAL = {Proceedings of the American Mathematical Society},
    VOLUME = {150},
      YEAR = {2022},
    NUMBER = {10},
     PAGES = {4325--4333},
      ISSN = {0002-9939,1088-6826},
   MRCLASS = {49J52 (30L15 35F30 37C10 49J35)},
  MRNUMBER = {4470177},
MRREVIEWER = {Jes\'us\ A.\ Jaramillo},
       DOI = {10.1090/proc/15958},
       URL = {https://doi.org/10.1090/proc/15958},
}

@article{LVV_2025,
  author    = {Felipe Lara and Nguyen Van Tuyen and Tran Van Nghi},
  title     = {Weak sharp minima at infinity and solution stability in mathematical programming via asymptotic analysis},
  journal   = {Journal of Global Optimization},
  year      = {2025},
  doi       = {10.1007/s10898-025-01516-2},
  url       = {https://doi.org/10.1007/s10898-025-01516-2}
}

@article {DLS_2024,
    AUTHOR = {Daniilidis, Aris and Le, Tri Minh and Salas, David},
     TITLE = {Metric compatibility and determination in complete metric
              spaces},
   JOURNAL = {Math. Z.},
  FJOURNAL = {Mathematische Zeitschrift},
    VOLUME = {308},
      YEAR = {2024},
    NUMBER = {4},
     PAGES = {Paper No. 62, 31},
      ISSN = {0025-5874,1432-1823},
   MRCLASS = {49J52 (03E10 30L15 49J53 58E05)},
  MRNUMBER = {4812490},
       DOI = {10.1007/s00209-024-03609-2},
       URL = {https://doi.org/10.1007/s00209-024-03609-2},
}

@article {E_1974,
    AUTHOR = {Ekeland, I.},
     TITLE = {On the variational principle},
   JOURNAL = {J. Math. Anal. Appl.},
  FJOURNAL = {Journal of Mathematical Analysis and Applications},
    VOLUME = {47},
      YEAR = {1974},
     PAGES = {324--353},
      ISSN = {0022-247X},
   MRCLASS = {49A25 (58E05)},
  MRNUMBER = {346619},
MRREVIEWER = {R.\ S.\ Palais},
       DOI = {10.1016/0022-247X(74)90025-0},
       URL = {https://doi.org/10.1016/0022-247X(74)90025-0},
}

@book {JJT_1983,
    AUTHOR = {Jongen, H. Th. and Jonker, P. and Twilt, F.},
     TITLE = {Nonlinear optimization in {${\bf R}^n$}. {I}},
    SERIES = {Methoden und Verfahren der Mathematischen Physik [Methods and
              Procedures in Mathematical Physics]},
    VOLUME = {29},
      NOTE = {Morse theory, Chebyshev approximation},
 PUBLISHER = {Verlag Peter D. Lang, Frankfurt am Main},
      YEAR = {1983},
     PAGES = {vi+263},
      ISBN = {3-8204-7903-1},
   MRCLASS = {90-02 (49M37 58C25 58E25 90C30)},
  MRNUMBER = {1072359},
MRREVIEWER = {W.\ Krabs},
}

@article {JW_1990,
    AUTHOR = {Jongen, Hubertus Th. and Weber, Gerhard-W.},
     TITLE = {On parametric nonlinear programming},
   JOURNAL = {Ann. Oper. Res.},
  FJOURNAL = {Annals of Operations Research},
    VOLUME = {27},
      YEAR = {1990},
    NUMBER = {1-4},
     PAGES = {253--283},
      ISSN = {0254-5330,1572-9338},
   MRCLASS = {90C31},
  MRNUMBER = {1088995},
       DOI = {10.1007/BF02055198},
       URL = {https://doi.org/10.1007/BF02055198},
}

@article {JW_1991,
    AUTHOR = {Jongen, H. Th. and Weber, G.-W.},
     TITLE = {Nonlinear optimization: characterization of structural
              stability},
   JOURNAL = {J. Global Optim.},
  FJOURNAL = {Journal of Global Optimization. An International Journal
              Dealing with Theoretical and Computational Aspects of Seeking
              Global Optima and Their Applications in Science, Management
              and Engineering},
    VOLUME = {1},
      YEAR = {1991},
    NUMBER = {1},
     PAGES = {47--64},
      ISSN = {0925-5001,1573-2916},
   MRCLASS = {90C31},
  MRNUMBER = {1263838},
MRREVIEWER = {Guang\ Jun\ Chen},
       DOI = {10.1007/BF00120665},
       URL = {https://doi.org/10.1007/BF00120665},
}

@article {JTW_1992,
    AUTHOR = {Jongen, H. T. and Twilt, F. and Weber, G.-W.},
     TITLE = {Semi-infinite optimization: structure and stability of the
              feasible set},
   JOURNAL = {J. Optim. Theory Appl.},
  FJOURNAL = {Journal of Optimization Theory and Applications},
    VOLUME = {72},
      YEAR = {1992},
    NUMBER = {3},
     PAGES = {529--552},
      ISSN = {0022-3239,1573-2878},
   MRCLASS = {90C34 (49M39 90C31)},
  MRNUMBER = {1147531},
MRREVIEWER = {Thomas\ Fischer},
       DOI = {10.1007/BF00939841},
       URL = {https://doi.org/10.1007/BF00939841},
}

@article {W_1998,
    AUTHOR = {Weber, Gerhard-W.},
     TITLE = {On the topology of parametric optimal control},
   JOURNAL = {J. Austral. Math. Soc. Ser. B},
  FJOURNAL = {Australian Mathematical Society. Journal. Series B. Applied
              Mathematics},
    VOLUME = {39},
      YEAR = {1998},
    NUMBER = {4},
     PAGES = {463--497},
      ISSN = {0334-2700},
   MRCLASS = {49K40 (90C31)},
  MRNUMBER = {1621329},
MRREVIEWER = {Hubertus\ Th.\ Jongen},
       DOI = {10.1017/S033427000000775X},
       URL = {https://doi.org/10.1017/S033427000000775X},
}

@article {PR_1998,
    AUTHOR = {Poliquin, R. A. and Rockafellar, R. T.},
     TITLE = {Tilt stability of a local minimum},
   JOURNAL = {SIAM J. Optim.},
  FJOURNAL = {SIAM Journal on Optimization},
    VOLUME = {8},
      YEAR = {1998},
    NUMBER = {2},
     PAGES = {287--299},
      ISSN = {1052-6234,1095-7189},
   MRCLASS = {90C31 (49J52 49K40)},
  MRNUMBER = {1618790},
MRREVIEWER = {B.\ Mordukhovich},
       DOI = {10.1137/S1052623496309296},
       URL = {https://doi.org/10.1137/S1052623496309296},
}

@article {HYZ_2024,
    AUTHOR = {Hu, Chunhai and Yang, Xiaoqi and Zheng, Xi Yin},
     TITLE = {Uniform weak sharp minima for multiobjective optimization
              problems},
   JOURNAL = {SIAM J. Optim.},
  FJOURNAL = {SIAM Journal on Optimization},
    VOLUME = {34},
      YEAR = {2024},
    NUMBER = {4},
     PAGES = {3699--3722},
      ISSN = {1052-6234,1095-7189},
   MRCLASS = {90C29 (49K40 90C31 90C48)},
  MRNUMBER = {4833223},
MRREVIEWER = {Xiang-Kai\ Sun},
       DOI = {10.1137/23M1628012},
       URL = {https://doi.org/10.1137/23M1628012},
}

@article {KMM_2019,
    AUTHOR = {Karkhaneei, Mohammad Mahdi and Mahdavi-Amiri, Nezam},
     TITLE = {Nonconvex weak sharp minima on {R}iemannian manifolds},
   JOURNAL = {J. Optim. Theory Appl.},
  FJOURNAL = {Journal of Optimization Theory and Applications},
    VOLUME = {183},
      YEAR = {2019},
    NUMBER = {1},
     PAGES = {85--104},
      ISSN = {0022-3239,1573-2878},
   MRCLASS = {49J52 (90C26 90C46)},
  MRNUMBER = {3989298},
MRREVIEWER = {Letizia\ Pellegrini},
       DOI = {10.1007/s10957-019-01539-2},
       URL = {https://doi.org/10.1007/s10957-019-01539-2},
}

@article {LMPHY_2018,
    AUTHOR = {Li, Chong and Meng, Li and Peng, Lihui and Hu, Yaohua and Yao,
              Jen-Chih},
     TITLE = {Weak sharp minima for convex infinite optimization problems in
              normed linear spaces},
   JOURNAL = {SIAM J. Optim.},
  FJOURNAL = {SIAM Journal on Optimization},
    VOLUME = {28},
      YEAR = {2018},
    NUMBER = {3},
     PAGES = {1999--2021},
      ISSN = {1052-6234,1095-7189},
   MRCLASS = {49K40 (46N10 49J52 90C25 90C31 90C34)},
  MRNUMBER = {3826677},
MRREVIEWER = {Adam\ B.\ Levy},
       DOI = {10.1137/16M1139564},
       URL = {https://doi.org/10.1137/16M1139564},
}

@article {X_2018,
    AUTHOR = {Xuan Duc Ha Truong},
     TITLE = {Slopes, error bounds and weak sharp {P}areto minima of a
              vector-valued map},
   JOURNAL = {J. Optim. Theory Appl.},
  FJOURNAL = {Journal of Optimization Theory and Applications},
    VOLUME = {176},
      YEAR = {2018},
    NUMBER = {3},
     PAGES = {634--649},
      ISSN = {0022-3239,1573-2878},
   MRCLASS = {49J53 (58C06 58E17 90C29)},
  MRNUMBER = {3772978},
MRREVIEWER = {Adam\ B.\ Levy},
       DOI = {10.1007/s10957-018-1240-6},
       URL = {https://doi.org/10.1007/s10957-018-1240-6},
}

@article {MOYZ_2025,
    AUTHOR = {Ma, Xiaoxiao and Ouyang, Wei and Ye, Jane J. and Zhang,
              Binbin},
     TITLE = {On second-order weak sharp minima of general nonconvex
              set-constrained optimization problems},
   JOURNAL = {J. Optim. Theory Appl.},
  FJOURNAL = {Journal of Optimization Theory and Applications},
    VOLUME = {207},
      YEAR = {2025},
    NUMBER = {2},
     PAGES = {Paper No. 21, 24},
      ISSN = {0022-3239,1573-2878},
   MRCLASS = {90C26 (49J52 49J53 90C46)},
  MRNUMBER = {4933996},
MRREVIEWER = {M.\ Alavi Hejazi},
       DOI = {10.1007/s10957-025-02775-5},
       URL = {https://doi.org/10.1007/s10957-025-02775-5},
}

@article {F_1991,
    AUTHOR = {Ferris, Michael C.},
     TITLE = {Finite termination of the proximal point algorithm},
   JOURNAL = {Math. Programming},
  FJOURNAL = {Mathematical Programming},
    VOLUME = {50},
      YEAR = {1991},
    NUMBER = {3},
     PAGES = {359--366},
      ISSN = {0025-5610,1436-4646},
   MRCLASS = {90C30 (49M37)},
  MRNUMBER = {1114237},
MRREVIEWER = {Stefan\ Mititelu},
       DOI = {10.1007/BF01594944},
       URL = {https://doi.org/10.1007/BF01594944},
}

@article {A_2006,
    AUTHOR = {Az\'e, D. and Corvellec, J.-N.},
     TITLE = {Variational methods in classical open mapping theorems},
   JOURNAL = {J. Convex Anal.},
  FJOURNAL = {Journal of Convex Analysis},
    VOLUME = {13},
      YEAR = {2006},
    NUMBER = {3-4},
     PAGES = {477--488},
      ISSN = {0944-6532,2363-6394},
   MRCLASS = {49J53 (46T20 54C30)},
  MRNUMBER = {2291548},
MRREVIEWER = {Russell\ Luke},
}
\bibliographystyle{plain}
\nocite{*}

\end{document}